# Network-Based Marketing: Identifying Likely Adopters via Consumer Networks

**Shawndra Hill, Foster Provost and Chris Volinsky**

*Abstract.* *Network-based marketing* refers to a collection of marketing techniques that take advantage of links between consumers to increase sales. We concentrate on the consumer networks formed using direct interactions (e.g., communications) between consumers. We survey the diverse literature on such marketing with an emphasis on the statistical methods used and the data to which these methods have been applied. We also provide a discussion of challenges and opportunities for this burgeoning research topic. Our survey highlights a gap in the literature. Because of inadequate data, prior studies have not been able to provide direct, statistical support for the hypothesis that network linkage can directly affect product/service adoption. Using a new data set that represents the adoption of a new telecommunications service, we show very strong support for the hypothesis. Specifically, we show three main results: (1) "Network neighbors"—those consumers linked to a prior customer—adopt the service at a rate 3–5 times greater than baseline groups selected by the best practices of the firm's marketing team. In addition, analyzing the network allows the firm to acquire new customers who otherwise would have fallen through the cracks, because they would not have been identified based on traditional attributes. (2) Statistical models, built with a very large amount of geographic, demographic and prior purchase data, are significantly and substantially improved by including network information. (3) More detailed network information allows the ranking of the network neighbors so as to permit the selection of small sets of individuals with very high probabilities of adoption.

## 1. INTRODUCTION

*Network-based marketing* seeks to increase brand recognition and profit by taking advantage of a so-cial network among consumers. Instances of network-based marketing have been called *word-of-mouth marketing*, *diffusion of innovation*, *buzz marketing* and *viral marketing* (we do not consider multilevel marketing, which has become known as "network" marketing). Awareness or adoption spreads from consumer to consumer. For example, friends or acquaintances may tell each other about a product or service, increasing awareness and possibly exercising explicit advocacy. Firms may use their websites to facilitate consumer-to-consumer advocacy via product recommendations (Kautz, Selman and Shah, 1997) or via on-line customer feedback mechanisms (Dellarocas, 2003). Consumer networks may also provide leverage to the advertising or marketing strategy of the firm. For example, in this paper we show how analysis of a consumer network improves targeted marketing.

*Shawndra Hill is a Doctoral Candidate and Foster Provost is Associate Professor, Department of Information, Operations and Management Sciences, Leonard N. Stern School of Business, New York University, New York, New York 10012-1126, USA e-mail: shill@stern.nyu.edu; fprovost@stern.nyu.edu. Chris Volinsky is Director, Statistics Research Department, AT&T Labs Research, Shannon Laboratory, Florham Park, New Jersey 07932, USA e-mail: volinsky@research.att.com.*

This paper makes two contributions. First we survey the burgeoning methodological research literature on network-based marketing, in particular on statistical analyses for network-based marketing. We review the research questions posed, and the data and analytic techniques used. We also discuss challenges and opportunities for research in this area. The review allows us to postulate necessary data requirements for studying the effectiveness of network-based marketing and to highlight the lack of current research that satisfies those requirements. Specifically, research must have access both to direct links between consumers and to direct information on the consumers' product adoption. Because of inadequate data, prior studies have not been able to provide direct, statistical support (Van den Bulte and Lilien, 2001) for the hypothesis that network linkage can directly affect product/service adoption.

The second contribution is to provide empirical support that network-based marketing indeed can improve on traditional marketing techniques. We introduce telecommunications data that present a natural testbed for network-based marketing models, in which communication linkages as well as product adoption rates can be observed. For these data, we show three main results: (1) "Network neighbors"—those consumers linked to a prior customer—adopt the service at a rate 3–5 times greater than baseline groups selected by the best practices of the firm's marketing team. In addition, analyzing the network allows the firm to acquire new customers who otherwise would have fallen through the cracks, because they would not have been identified based on traditional attributes. (2) Statistical models, built with a very large amount of geographic, demographic and prior purchase data, are significantly and substantially improved by including network information. (3) More sophisticated network information allows the ranking of the network neighbors so as to permit the selection of small sets of individuals with very high probabilities of adoption.

## 2. NETWORK-BASED MARKETING

There are three, possibly complementary, modes of network-based marketing.

*Explicit advocacy*: Individuals become vocal advocates for the product or service, recommending it to their friends or acquaintances. Particular individuals such as Oprah, with her monthly book club reading list, may represent "hubs" of advocacy in the consumer relationship network. The success of *The Da Vinci Code*, by Dan Brown, may be due to its initial marketing: 10,000 books were delivered free to readers thought to be influential enough (e.g., individuals, booksellers) to stimulate the traffic in paid-for editions (Paumgarten, 2003). When firms give explicit incentives to consumers to spread information about a product via word of mouth, it has been called *viral marketing*, although that term could be used to describe any network-based marketing where the pattern of awareness or adoption spreads from consumer to consumer.

*Implicit advocacy*: Even if individuals do not speak about a product, they may advocate implicitly through their actions—especially through their own adoption of the product. Designer labeling has a long tradition of using consumers as implicit advocates. Firms commonly capitalize on influential individuals (such as athletes) to advocate products simply by conspicuous adoption. More recently, firms have tried to induce the same effect by convincing particularly "cool" members of smaller social groups to adopt products (Gladwell, 1997; Hightower, Brady and Baker, 2002).

*Network targeting*: The third mode of network-based marketing is for the firm to market to prior purchasers' social-network neighbors, possibly without any advocacy at all by customers. For network targeting, the firm must have some means to identify these social neighbors.

These three modes may be used in combination. A well-cited example of viral marketing combines network targeting and implicit advocacy: The Hotmail free e-mail service appended to the bottom of every outgoing e-mail message the hyperlinked advertisement, "Get your free e-mail at Hotmail," thereby targeting the social neighbors of every current user (Montgomery, 2001), while taking advantage of the user's implicit advocacy. Hotmail saw an exponentially increasing customer base. Started in July 1996, in the first month alone Hotmail acquired 20,000 customers. By September 1996 the firm had acquired over 100,000 accounts, and by early 1997 it had over 1 million subscribers.

Traditional marketing methods do not appeal to some segments of consumers. Some consumers apparently value the appearance of being on the cutting edge or "in the know," and therefore derive satisfaction from promoting new, exciting products. The firm BzzAgents (Walker, 2004) has managed to entice voluntary (unpaid) marketing of new products. Furthermore, although more and more information has become available on products, parsing such information is costly to the consumer. Explicit advocacy, such as word-of-mouth advocacy, can be a useful way to filter out noise.

A key assumption of network-based marketing through explicit advocacy is that consumers propagate "positive" information about products after they either have been made aware of the product by traditional marketing vehicles or have experienced the product themselves. Under this assumption, a particular subset of consumers may have greater value to firms because they have a higher propensity to propagate product information (Gladwell, 2002), based on a combination of their being particularly influential and their having more friends (Richardson and Domingos, 2002). Firms should want to find these influencers and to promote useful behavior.

## 3. LITERATURE REVIEW

Many quantitative statistical methods used in empirical marketing research assume that consumers act independently. Typically, many explanatory attributes are collected on each actor and used in multivariate modeling such as regression or tree induction. In contrast, network-based marketing assumes interdependency among consumer preferences. When interdependencies exist, it may be beneficial to account for their effects in targeting models. Traditionally in statistical research, interdependencies are modeled as part of a covariance structure, either within a particular observational unit (as in the case of repeated measures experiments) or between observational units. Studies of network-based marketing instead attempt to measure these interdependencies through implicit links, such as matching on geographic or demographic attributes, or through explicit links, such as direct observation of communications between actors. In this section, we review the different types of data and the range of statistical methods that have been used to analyze them, and we discuss the extent to which these methods naturally accommodate networked data.

Work in network-based marketing spans the fields of statistics, economics, computer science, sociology, psychology and marketing. In this section, we organize prominent work in network-based marketing by six types of statistical research: (1) econometric modeling, (2) network classification modeling, (3) surveys, (4) designed experiments with convenience samples, (5) diffusion theory and (6) collaborative filtering and recommender systems. In each case, we provide an overview of the approach and a discussion of a prominent example. This (brief) survey is not exhaustive. In the final subsection, we discuss some of the statistical challenges inherent in incorporating this network structure.

### 3.1 Econometric Models

Econometrics is the application of statistical methods to the empirical estimation of economic relationships. In marketing this often means the estimation of two simultaneous equations: one for the marketing organization or firm and one for the market. Regression and time-series analysis are found at the core of econometric modeling, and econometric models are often used to assess the impact of a target marketing campaign over time.

Econometric models have been used to study the impact of interdependent preferences on rice consumption (Case, 1991), automobile purchases (Yang and Allenby, 2003) and elections (Linden, Smith and York, 2003). For each of the aforementioned studies, geography is used in part as a proxy for interdependence between consumers, as opposed to direct, explicit communication. However, different methods are used in the analysis. Most recently, Yang and Allenby (2003) suggested that traditional random effects models are not sufficient to measure the interdependencies of consumer networks. They developed a Bayesian hierarchical mixture model where interdependence is built into the covariance structure through an autoregressive process. This framework allows testing of the presence of interdependence through a single parameter. It also can incorporate the effects of multiple networks, each with its own estimated dependence structure. In their application, they use geography and demography to create a "network" of consumers in which links are created between consumers who exhibit geographic or demographic similarity. The authors showed that the geographically defined network of consumers is more useful than the demographic network for explaining consumer behavior as it relates to purchasing Japanese cars. Although they do not have data on direct communication between consumers, the framework presented by Yang and Allenby (2003) could be extended to explicit network data where links are created between consumers through their explicit communication as opposed to demographic or geographic similarity.

A drawback of this approach is that the interdependence matrix has size $n^2$, where $n$ is the number of consumers; consumer networks are extremely large and prohibit parameter estimation using this method. Sparse matrix techniques or clever clustering of the observations would be a natural extension.

### 3.2 Network Classification Models

Network classification models use knowledge of the links between entities in a network to estimate a quantity of interest for those entities. Typically, in such a model an entity is influenced most by those directly connected to it, but is also affected to a lesser extent by those further away. Some network classification models use an entire network to make predictions about a particular entity on the network; Macskassy and Provost (2004) provided a brief survey. However, most methods have been applied to small data sets and have not been applied to consumer data. Much research in network classification has grown out of the pioneering work by Kleinberg (1999) on hubs and authorities on the Internet, and out of Google's PageRank algorithm (Brin and Page, 1998), which (to oversimplify) identifies the most influential members of a network by how many influential others "point" to them. Although neither study uses statistical models, both are related to well-understood notions of degree centrality and distance centrality from the field of social-network analysis.

One paper that models a consumer network for maximizing profit is by Richardson and Domingos (2002), in which a social network of customers is modeled as a Markov random field. The probability that a given customer will buy a given product is a function of the states of her neighbors, attributes of the product and whether or not the customer was marketed to. In this framework it is possible to assign a "network value" to every customer by estimating the overall benefit of marketing to that customer, including the impact that the marketing action will have on the rest of the network (e.g., through word of mouth). The authors tested their model on a database of movie reviews from an Internet site and found that their proposed methodology outperforms non-network methods for estimating customer value. Their network formulation uses implicit links (customers are linked when a customer reads a review by another customer and subsequently reviews the item herself) and implicit purchase information (they assume a review of an item implies a purchase and vice versa).

### 3.3 Surveys

Most research in this area does not have information on whether consumers actually talk to each other. To address this shortcoming, some studies use survey sampling to collect comprehensive data on consumers' word-of-mouth behavior. By sampling individuals and contacting them, researchers can collect data that are difficult (or impossible) to obtain directly by observing network-based marketing phenomena (Bowman and Narayandas, 2001). The strength of these studies lies in the data, including the richness and flexibility of the answers that can be collected from the responders. For instance, researchers can acquire data about how customers found out about a product and how many others they told about the product. An advantage is that researchers can design their sampling scheme to control for any known confounding factors and can devise fully balanced experimental designs that test their hypotheses. Since the purpose of models built from survey data is description, simple statistical methods like logistic regression or analysis of variance (ANOVA) typically are used.

Bowman and Narayandas (2001) surveyed more than 1700 purchasers of 60 different products who previously had contacted the manufacturer of that product. The purchasers were asked specific questions about their interaction with the manufacturer and its impact on subsequent word-of-mouth behavior. The authors were able to capture whether the customers told others of their experience and if so, how many people they told. The authors found that self-reported "loyal" customers were more likely to talk to others about the products when they were dissatisfied, but interestingly not more likely when they were satisfied. Although studies like this collect some direct data on consumers' word-of-mouth behavior, the researchers do not know which of the consumers' contacts later purchased the product. Therefore, they cannot address whether word-of-mouth actually affects individual sales.

### 3.4 Designed Experiments with Convenience Samples

Designed experiments enable researchers to study network-based marketing in a controlled setting. Although the subjects typically comprise a convenience sample (such as those undergraduates who answer an ad in the school newspaper), the design of the experiment can be completely randomized. This is unlike the studies that rely on secondary data sources or data from the Web. Typically ANOVA is used to draw conclusions.

Frenzen and Nakamoto (1993) studied the factors that influence individuals' decisions to disseminate information through a market via word-of-mouth. The subjects were presented with several scenarios that represented different products and marketing strategies, and were asked whether they would tell trusted and nontrusted acquaintances about the product/sale. They studied the effect of the cost/value manipulations on the consumers' willingness to share information actively with others, as a function of the strength of the social tie. In this study, the authors did not allow the subjects to construct their explicit consumer network; instead, they asked the participants to hypothesize about their networks. The experiments used the data from a convenience sample to generalize over a complete consumer network. The authors also employed simulations in their study. They found that the stronger the moral hazard (the risk of problematic behavior) presented by the information, the stronger the ties must be to foster information propagation. Generally, the authors showed that network structure and information characteristics interact when individuals form their information transmission decisions.

### 3.5 Diffusion Models

Diffusion theory provides tools, both quantitative and qualitative, to assess the likely rate of diffusion of a technology or product. Qualitatively, researchers have identified numerous factors that facilitate or hinder technology adoption (Fichman, 2004), as well as social factors that influence product adoption (Rogers, 2003). Quantitative diffusion research involves empirical testing of predictions from diffusion models, often informed by economic theory.

The most notable and most influential diffusion model was proposed by Bass (1969). The Bass model of product diffusion predicts the number of users who will adopt an innovation at a given time $t$. It hypothesizes that the rate of adoption is a function solely of the current proportion of the population who have adopted. Specifically, let $F(t)$ be the cumulative proportion of adopters in the population. The diffusion equation, in its simplest form, models $F(t)$ as a function of $p$, the intrinsic adoption rate, and $q$, a measure of social contagion. When $q > p$, this equation describes an $S$-shaped curve, where adoption is slow at first, takes off exponentially and tails off at the end. This model can effectively model word-of-mouth product diffusion at the aggregate, societal level.

In general, the empirical studies that test and extend accepted theories of product diffusion rely on aggregate-level data for both the customer attributes and the overall adoption of the product (Ueda, 1990; Tout, Evans and Yakan, 2005); they typically do not incorporate individual adoption. Models of product diffusion assume that network-based marketing is effective. Since understanding when diffusion occurs and the extent to which it is effective is important for marketers, these methods may benefit from using individual-level data. Data on explicit networks would enable the extension of existing diffusion models, as well as the comparison of results using individual- versus aggregate-level data.

In his first study, Bass tested his model empirically against data for 11 consumer durables. The model yielded good predictions of the sales peak and the timing of the peak when applied to historical data. Bass used linear regression to estimate the parameters for future sales predictions, measuring the goodness of fit ($R^2$ value) of the model for 11 consumer durable products. The success of the forecasts suggests that the model may be useful in providing long-range forecasting for product sales or adoption. There has been considerable follow-up work on diffusion since this groundbreaking work. Mahajan, Muller and Kerin (1984) review this work. Recent work on product diffusion explores the extent to which the Internet (Fildes, 2003) as well as globalization (Kumar and Krishnan, 2002) play a role in product diffusion.

### 3.6 Collaborative Filtering and Recommender Systems

Recommender systems make personalized recommendations to individual consumers based on demographic content and link data (Adomavicius and Tuzhilin, 2005). Collaborative filtering methods focus on the links between consumers; however, the links are not direct. They associate consumers with each other based on shared purchases or similar ratings of shared products.

Collaborative filtering is related to explicit consumer network-based marketing because both target marketing tasks benefit from learning from data stored in multiple tables (Getoor, 2005). For example, Getoor and Sahami (1999), Huang, Chung and Chen (2004) and Newton and Greiner (2004) established the connection between the recommendation problem and statistical relational learning through the application of probabilistic relational models (PRM's) (Getoor, Friedman, Koller and Pfeffer, 2001). However, neither group used explicit links between customers for learning. Recommendation systems may well benefit from information about explicit consumer interaction as an additional, perhaps quite important, aspect of similarity.

### 3.7 Research Opportunities and Statistical Challenges

We see that there is a burgeoning body of work that addresses consumers' interactions and their effects on purchasing. To our knowledge the foregoing types represent the main statistical approaches taken in research on network-based marketing. In each approach, there are assumptions made in the data collection or in the analysis that restrict them from providing strong and direct support for the hypothesis that network-based marketing indeed can improve on traditional techniques. Surveys and convenience samples can suffer from small and possibly biased samples. Collaborative filtering models have large samples, but do not measure direct links between individuals. Models in network classification and econometrics historically have used proxies like geography instead of data on direct communications, and almost all studies have no accurate, specific data on which (and what) customers purchase.

To paint a complete picture of network influence for a particular product, the ideal data set would have the following properties: (1) large and unbiased sample, (2) comprehensive covariate information on subjects, (3) measurement of direct communication between subjects and (4) accurate information on subjects' purchases. The data set we present in the next section has all of these properties and we will demonstrate its value for statistical research into network influence. The question of how to analyze such data brings up many statistical issues:

*Data-set size.* Network-based marketing data sets often arise from Internet or telecommunications applications and can be quite large. When observations number in the millions (or hundreds of millions), the data become unwieldy for the typical data analyst and often cannot be handled in memory by standard statistical analysis software. Even if the data can be loaded, their size renders the interactive style of analysis common with tools like R or Splus painfully slow. In Internet or telecommunications studies, there often are two massive sources of data: all actors (web sites, communicators), along with their descriptive attributes, and the transactions among these actors. One solution is to compress the transaction information into attributes to be included in the actors' attribute set. It has been shown that file squashing (DuMouchel, Volinsky et al., 1999), which attempts to combine the best features of preprocessed data with random sampling, can be useful for customer attrition prediction. DuMouchel et al. claimed that squashing can be useful when dealing with up to billions of records. However, there may be a loss of important information which can be captured only by complex network structure.

More sophisticated network information derived from transactional data can also be incorporated into the matrix of customer information by deriving network attributes such as degree distribution and time spent on the network (which we demonstrate below). Similarly, other types of data such as geographical data or temporal data, which otherwise would need to be handled by some sophisticated methodology, can be folded into the analysis by creating new covariates. It remains an open question whether clever data engineering can extract all useful information to create a set of covariates for traditional analysis. For example, knowledge of communication with specific sets of individuals can be incorporated, and may provide substantial benefit (Perlich and Provost, 2006).

Once the data are combined, the remaining data set still may be quite large. While much data mining research is focused on scaling up the statistical toolbox to today's massive data sets, random sampling remains an effective way to reduce data to a manageable size while maintaining the relationships we are trying to discover, if we assume the network information is fully encoded in the derived variables. The amount of sampling necessary will depend on the computing environment and the complexity of the model, but most modern systems can handle data sets of tens or hundreds of thousands of observations. When sampling, care must be taken to stratify by any attributes that are of particular interest or to oversample those attributes that have extremely skewed distributions.

*Low incidence of response.* In applications where the response is a consumer's purchase or reaction to a marketing event, it is common to have a very low response rate, which can result in poor fit and reduced ability to detect significant effects for standard techniques like logistic regression. If there are not many independent attributes, one solution is Poisson regression, which is well suited for rare events. Poisson regression requires forming buckets of observations based on the independent attributes and modeling the aggregate response in these buckets as a Poisson random variable. This requires discretization of any continuous independent attributes, which may not be desirable. Also, if there are even a moderate number of independent attributes, the buckets will be too sparse to allow Poisson modeling. Other solutions that have been proposed include oversampling positive responses and/or undersampling negative responses. Weiss (2004) gave an overview of the literature on these and related techniques, showing that there is mixed evidence as to their effectiveness. Other studies of note include the following. Weiss and Provost (2003) showed that, given a fixed sample size, the optimal class proportion in training data varies by domain and by ultimate objective (but can be determined); generally speaking, to produce probability estimates or rankings, a 50:50 distribution is a good default. However, Weiss and Provost's results are only for tree induction. Japkowicz and Stephen (2002) experimented with neural networks and support-vector machines, in addition to tree induction, showing (among other things) that support-vector machines are insensitive to class imbalance. However, they considered primarily noise-free data. Other techniques to deal with unbalanced response attributes include ensemble (Chan and Stolfo, 1998; Mease, Wyner and Buja, 2006) and multiphase rule induction (Clearwater and Stern, 1991; Joshi, Kumar and Agarwal, 2001). This is an area in need of more systematic empirical and theoretical study.

*Separating word-of-mouth from homophily.* Unless there is information about the content of communications, one cannot conclude that there was word-of-mouth transmission of information about the product. Social theory tells us that people who communicate with each other are more likely to be similar to each other, a concept called *homophily* (Blau, 1977; McPherson, Smith-Lovin and Cook, 2001). Homophily is exhibited for a wide variety of relationships and dimensions of similarity. Therefore, linked consumers probably are like-minded, and like-minded consumers tend to buy the same products. One way to address this issue in the analysis is to account for consumer similarity using *propensity scores* (Rosenbaum and Rubin, 1984). Propensity scores were developed in the context of nonrandomized clinical trials and attempt to adjust for the fact that the statistical profile of patients who received treatment may be different than the profile of those who did not, and that these differences could mask or enhance the apparent effect of the treatment. Let $T$ represent the treatment, $\mathbf{X}$ represent the independent attributes excluding the treatment and $Y$ represent the response. Then the propensity score $\mathrm{PS}(\mathbf{x}) = P(T = 1|\mathbf{X} = \mathbf{x})$. By matching propensity scores in the treatment and control groups using typical indicators of homophily like demographic data, we can account (partially) for the possible confoundedness of other independent attributes.

*Incorporating extended network structure.* Data with network structure lend themselves to a robust set of network-centric analyses. One simple method (employed in our analysis) is to create attributes from the network data and plug them into a traditional analysis. Another approach is to let each actor be influenced by her neighborhood modeled as a Markov random field. Domingos and Richardson (2001) used this technique to assign every node a "network value." A node with high network value (1) has a high probability of purchase, (2) is likely to give the product a high rating, (3) is influential on its neighbors' ratings and (4) has neighbors like itself. Hoff, Raftery and Handcock (2002) defined a Markov-chain Monte Carlo method to estimate latent positions of the actors for small social-network data sets. This embeds the actors in an unobserved "social space," which could be more useful than the actual transactions themselves for predicting sales. The field of statistical relational learning (Getoor, 2005) has recently produced a wide variety of methods that could be applicable. Often these models allow influence to propagate through the network.

*Missing data.* Missing data in network transactions are common—often only part of a network is observable. For instance, firms typically have transactional data on their customers only or may have one class of communication (e-mail) but not another (cellular phone). One attempt to account for these missing edges is to use network structure to assign a probability of a missing edge everywhere an edge is not present. Thresholding this probability creates *pseudo-edges*, which can be added to the network, perhaps with a lesser weight (Agarwal and Pregibon, 2004). This is closely related to the *link prediction* problem, which tries to predict where the next links will be (Liben-Nowell and Kleinberg, 2003). One extension of the PRM framework models link structure through the use of *reference uncertainty* and *existence uncertainty*. The extension includes a unified generative model for both content and relational structure, where interactions between the attributes and link structure are modeled (Getoor, Friedman, Koller and Taskar, 2003).

## 4. DATA SET AND PRIMARY HYPOTHESIS

This section details our data set, derived primarily from a direct-mail marketing campaign to potential customers of a new communications service (later we augment the primary data with a large set of consumer-specific attributes). The firm's marketing team identified and marketed to a list of prospects using its standard methods. We investigate whether network-related effects or evidence of "viral" information spread are present in this group. As we will describe, the firm also marketed to a group we identified using the network data, which allows us to test our hypotheses further. We are not permitted to disclose certain details, including specifics about the service being offered and the exact size of the data set.

### 4.1 Initial Data Details

In late 2004, a telecommunications firm undertook a large direct-mail marketing campaign to potential customers of a new communications service. This service involved new technology and, because of this, it was believed that marketing would be most successful to those consumers who were thought to be "high tech."

TABLE 1
*Descriptive statistics for the marketing segments (see Section 4.1 for details)*

| Segment | Loyalty | Intl | Tech1 | Tech2 | Early Adopt | Offer | % of list | %NN |
|---|---|---|---|---|---|---|---|---|
| 1 | 3 | Y | Hi | 1–7 | Med–Hi | P1 | 1.6 | 0.63 |
| 2 | 3 | Y | Med | 1–7 | Med–Hi | P1 | 2.4 | 1.26 |
| 3 | 2 | Y | Hi | 1–4 | Hi | P1 | 1.7 | 0.08 |
| 4 | 2 | Y | Med | 1–4 | Hi | P1 | 1.7 | 0.10 |
| 5 | 1 | Y | Hi | 1–4 | Hi | P1 | 0.1 | 0.22 |
| 6 | 1 | Y | Med | 1–4 | Hi | P1 | 0.1 | 0.25 |
| 7 | 3 | N | Hi | 1–7 | Med–Hi | P2 | 10.9 | 0.50 |
| 8 | 3 | N | Med | 1–7 | Med–Hi | P2 | 13.1 | 0.83 |
| 9 | 2 | N | Hi | 1–4 | Hi | P2 | 17.5 | 0.04 |
| 10 | 2 | N | Med | 1–4 | Hi | P2 | 11.0 | 0.07 |
| 11 | 1 | N | Hi | 1–4 | Hi | P2 | 5.3 | 0.14 |
| 12 | 1 | N | Med | 1–4 | Hi | P2 | 7.7 | 0.25 |
| 13 | 3 | N | Hi | 1–7 | Med–Hi | P2 | 2.0 | 0.63 |
| 14 | 1, 2 | N | Hi | 1–4 | Hi | P2 | 2.0 | 0.15 |
| 15 | 1 | Y | ? | ? | ? | P3 | 2.0 | 1.01 |
| 16 | 1 | N | ? | ? | ? | P2 | 1.6 | 0.46 |
| 17 | 3 | N | Hi | 1–7 | Med–Hi | P2+ | 2.0 | 0.70 |
| 18 | 1, 2 | N | Hi | 1–4 | Hi | P2+ | 2.0 | 0.15 |
| 19 | 1, 2, 3 | Y | Hi | 1–7 | Med–Hi | P3 | 1.8 | 0.67 |
| 20 | 2 | N | Hi, Med | 1–4 | Hi | L1 | 6.0 | 0.05 |
| 21 | 2 | N | Hi, Med | 1–4 | Hi | L2 | 6.0 | 0.05 |

In keeping with standard practice, the marketing team collected attributes on a large set of prospects—consumers whom they believed to be potential adopters of the service. The marketing team used demographic data, customer relationship data, and various other data sources to create profitability and behavioral models to identify prospective *targets*—consumers who would receive a targeted mailing. The data the marketing team provided us with did not contain the underlying customer attributes (e.g., demographics), but instead included values for derived attributes that defined 21 marketing segments (Table 1) that were used for campaign management and post hoc analyses. The sample included millions of consumers. The team believed that the different segments would have varying response rates and it was important to separate the segments in this way to learn the most from the campaign.

An important derived variable was loyalty, a three-level score based on previous relationships with the firm, including previous orders of this and other services. Roughly, loyalty level 3 comprises customers with moderate-to-long tenure and/or those who have subscribed to a number of services in the past. Loyalty level 2 comprises those customers with which the firm has had some limited prior experiences. Loyalty level 1 comprises consumers who did not have service with the firm at the time of mailing; little (if any) information is available on them. Previous analyses have shown that loyalty and tenure attributes have substantial impact on response to campaigns.

Other important attributes were based on demographics and other customer characteristics. The attribute Intl is an indicator of whether the prospect had previously ordered any international services; Tech1 (hi, med or low) and Tech2 (1–10, where 1 = high tech) are scores derived from demographics and other attributes that estimate the interest and ability of the customer to use a high-tech service; Early Adopt is a proprietary score that estimates the likelihood of the customer to use a new product, based on previous behavior. We also show the Offer, indicating that different segments received different marketing messages: P1–P3 indicate different postcards that were sent, L1 and L2 indicate different letters, and a "+" indicates that a "call blast" accompanied the mailing. In defining the segments, those groups with high loyalty values were permitted lower values from the technology and early adoption models. Segments 15 and 16 were provided by an external vendor; there were insufficient data on these prospects to fit our Tech and Early Adopt models, as indicated by a "?" in Table 1.

### 4.2 Primary Hypothesis and Network Neighbors

The research goal we consider here is whether relaxing the assumption of independence between consumers can improve demonstrably the estimation of response likelihood. Thus, our first hypothesis is that someone who has direct communication with a current subscriber is more likely herself to adopt the service. It should be noted that the firm knows only of communications initiated by one of its customers through a service of the firm, so the network data are incomplete (considerably), especially for the lower loyalty groups. Data on communications events include anonymous identifiers for the transactors, a time stamp and the transaction duration. For the purposes of this research, all data are rendered anonymous so that individual identities are protected.

In pursuit of our hypothesis, we constructed an attribute called *network neighbor* (or NN)—a flag that indicates whether the targeted consumer had communicated with a current user of the service in a time period prior to the marketing campaign. Overall, 0.3% of the targets are network neighbors. In Table 1, the percentage of network neighbors (%NN) is broken down by segment.

In addition, the marketing team invited us to create our own segment, which they also would target. Our "segment 22" consisted of network neighbors that were not already on the current list of targets. To make sure our list contained viable prospects, the marketing team calculated the derived technology and early adopter scores for the consumers on our list. They filtered based on these scores, but they relaxed the thresholds used to limit their original list. For instance, someone with loyalty = 1 needed a Tech2 score less than 4 to merit inclusion on the initial list; this threshold was relaxed for our list to Tech2 less than 7. In this way, the marketing team allowed prospects who missed inclusion on the first cut to make it into segment 22 if they were network neighbors. However, the marketing team still avoided targeting customers who they believed had

TABLE 2
*Data categories*

| | **Target = Y** | **Target = N** |
|---|---|---|
| NN = Y | *NN targets*<br>Segments 1–22<br>Relative size = 0.015<br>Prospects identified by marketing models and who also are network neighbors. Those in segment 22 have reduced thresholds on the marketing model scores. | *NN nontargets*<br><br>Relative size = 0.10<br>Consumers who were network neighbors, but were not marketed to because they scored poorly on marketing models. |
| NN = N | *Non-NN targets*<br>Segments 1–21<br>Relative size = 1<br>Prospects identified by marketing models but who are not network neighbors. | *Non-NN nontargets*<br><br>Relative size > 8<br>Consumers who were not network neighbors and also were not considered to be good prospects by the marketing model. |

NOTES. The data for our study are broken down into targets and network neighbors. The "relative size" value shows the number of prospects who show up in each group, relative to the non-NN target group.

very small probabilities of a purchase. For those network neighbors who did not score high enough to warrant inclusion in segment 22, we still tracked their purchase records to see if any of them subscribed to the service in the absence of the marketing campaign; see below. Overall, the profile of the candidates in our segment 22 was considered to be subpar in terms of demographics, affinity and technological capability. Notably, for our final conclusions, these targets are potential customers the firm would have otherwise ignored. The size of segment 22 was about 1.2% of the marketing list.

To summarize, the above process divides the prospect universe along two dimensions: (1) *targets*—those consumers identified by the marketing models as being worthy of solicitation—and (2) *network neighbors*—those who had direct communication with a subscriber. Table 2 shows the relative size for each combination (using the non-network-neighbor targets as the reference set). Note the non-NN nontargets, who neither are network neighbors nor are they deemed to be good prospects. This group is the majority of the prospect space and includes consumers that the firm has very little information about, because they are low-usage communicators or do not subscribe to any services with the firm.

## 4.3 Modeling with Consumer-Specific Data

To determine whether relaxing the independence assumption (using the network data) improves modeling, we fit models using a wide range of demographic and consumer-specific independent attributes (many of which are known or believed to affect the estimated likelihood of purchase). Overall, we collected the values for over 150 attributes to assess their effect on sales likelihood and their interactions with the network-neighbor variable. These values included the following:

- Loyalty data: We obtained finer-grained loyalty information than the simple categorization described above, including past spending, types of service, how often the customer responded to prior mailings, a loyalty score generated by a proprietary model and information about length of tenure.
- Geographic data: Geographic data were necessary for the direct mail campaign. These data include city, state, zip code, area code and metropolitan city code.
- Demographic data: These include information such as gender, education level, credit score, head of household, number of children in the household, age of members in the household, occupation and home ownership. Some of this information was inferred at the census tract level from the geographic data.
- Network attributes: As mentioned earlier, we observed communications of current subscribers with other consumers. In addition to the simple network-neighbor flag described earlier, we derived more sophisticated attributes from prospects' communication patterns. We will return to these in Section 5.6.

### 4.4 Data Limitations

We encountered missing values for customers across all loyalty levels. The amount of missing information is directly related to the level of experience we have had with the customer just prior to the direct mailing. For example, geography data are available for all targets across all three loyalty levels. On the other hand, as the number of services and tenure with the firm decline, so does the amount of information (e.g., transactions) available for each target. Given the difference in information as loyalty varies, we grouped customers by loyalty level and treated the levels separately in our analyses. This stratification leaves three groups that are mostly internally consistent with respect to missing values.

The overall response rate is very low. As discussed above, this presents challenges inherent with a heavily skewed response variable. For example, an analysis that stratifies over many different attributes may have several strata with no sales at all, rendering these strata mostly useless. The data set is large, which helps to ameliorate this problem, but in turn presents logistical problems with many sophisticated statistical analyses. In this paper, we restrict ourselves to relatively straightforward analyses.

### 4.5 Loyalty Distribution

A look at the distribution of the loyalty groups across the four categories (Figure 1) of prospects shows that the firm targeted customers in the higher loyalty groups relatively heavily. The network-neighbor target group appears to skew toward the less loyal prospects; this is due to the fact that segment 22, which makes up a large part of the network-neighbor population, comprises predominantly low-loyalty consumers.

## 5. ANALYSIS

Next we will show direct, statistical evidence that consumers who have communicated with prior customers are more likely to become customers. We show this in several ways, including using our own best efforts to build competing targeting models and conducting thorough assessments of predictive ability on out-of-sample data. Then we consider more sophisticated network attributes and show that targeting can be improved further.

### 5.1 Network-Based Marketing Improves Response

Segmentation provides an ideal setting to test the significance and magnitude of any improvement in modeling by including network-neighbor information, while stratifying by many attributes known to be important, such as loyalty and tenure. The response variable is the take rate for the targets in the two months following the direct mailing. The take rate is the proportion of the targeted consumers who adopted the service within a specified period following the offer. For each segment, we performed a simple logistic regression for the independent network-neighbor attribute versus the dependent sales response. In Figure 2, we graphically present parameter estimates (equivalent to log-odds ratios) for the network attribute along with 95% confidence intervals for 20 of the 21 segments (segment 5 had only a small number of network-neighbor prospects and zero network-neighbor sales, and therefore had an infinite log odds). Figure 2 shows that in all 20 segments the network-neighbor effect is positive (the parameter estimate is greater than zero), demonstrating an increased take rate for the network-neighbor group within each segment. For 17 of these segments, the log-odds ratio is significantly different from the null hypothesis value of 0 ($p < 0.05$), indicating that being a network neighbor significantly affected sales in those segments.

While odds ratios allow for tests of significance of an independent variable, they are not as directly interpretable as comparisons of take rates of the network-neighbor and non-network-neighbor groups in a given segment. The take rates for the network neighbors are plotted versus the non-network neighbors in Figure 3, where the size of the point is proportional to the log size of the segment. All segments have higher take rates in the network-neighbor subgroup, except for the one segment that had no network-neighbor sales (the smallest sample size). Over the entire data set, the network-neighbors' take rates were greater by a factor of 3.4. This value is plotted in Figure 3 as a dotted line with slope = 3.4. The right-hand plot of Figure 3 shows the relationship between each segment's take rate and its lift ratio, defined as the take rate for NN divided by the take rate for non-NN. The plot shows that the benefit of being a network neighbor is greater for those segments with lower overall take rates.

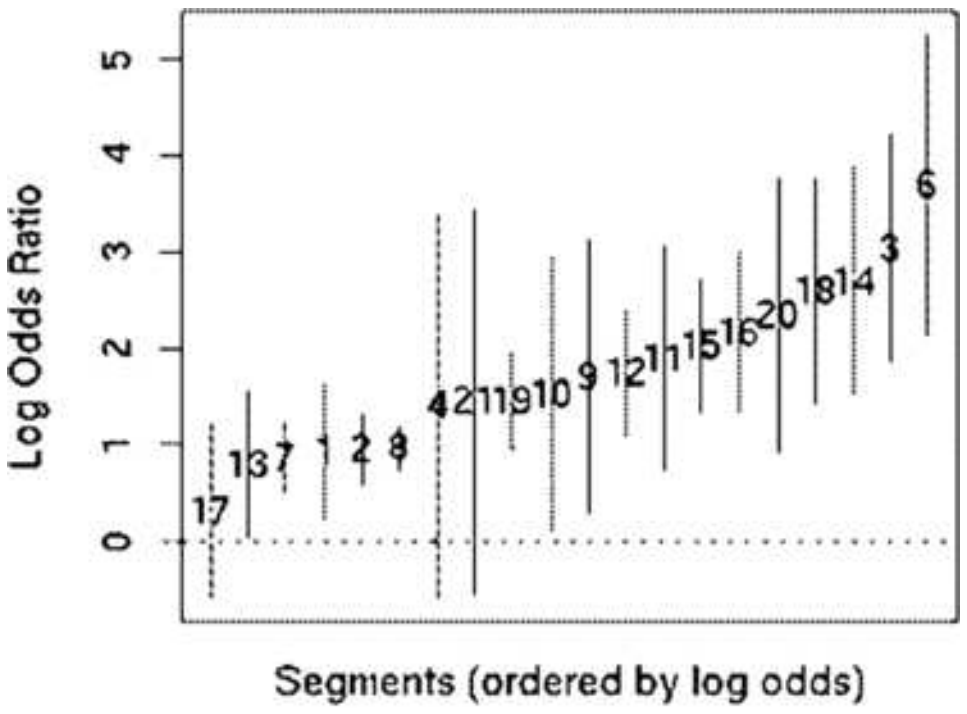

Fig. 2. *Results of logistic regression. Parameter estimates plotted as log-odds ratios with 95% confidence intervals. The number plotted at the value of the parameter estimate refers back to segment numbers from Table 1.*

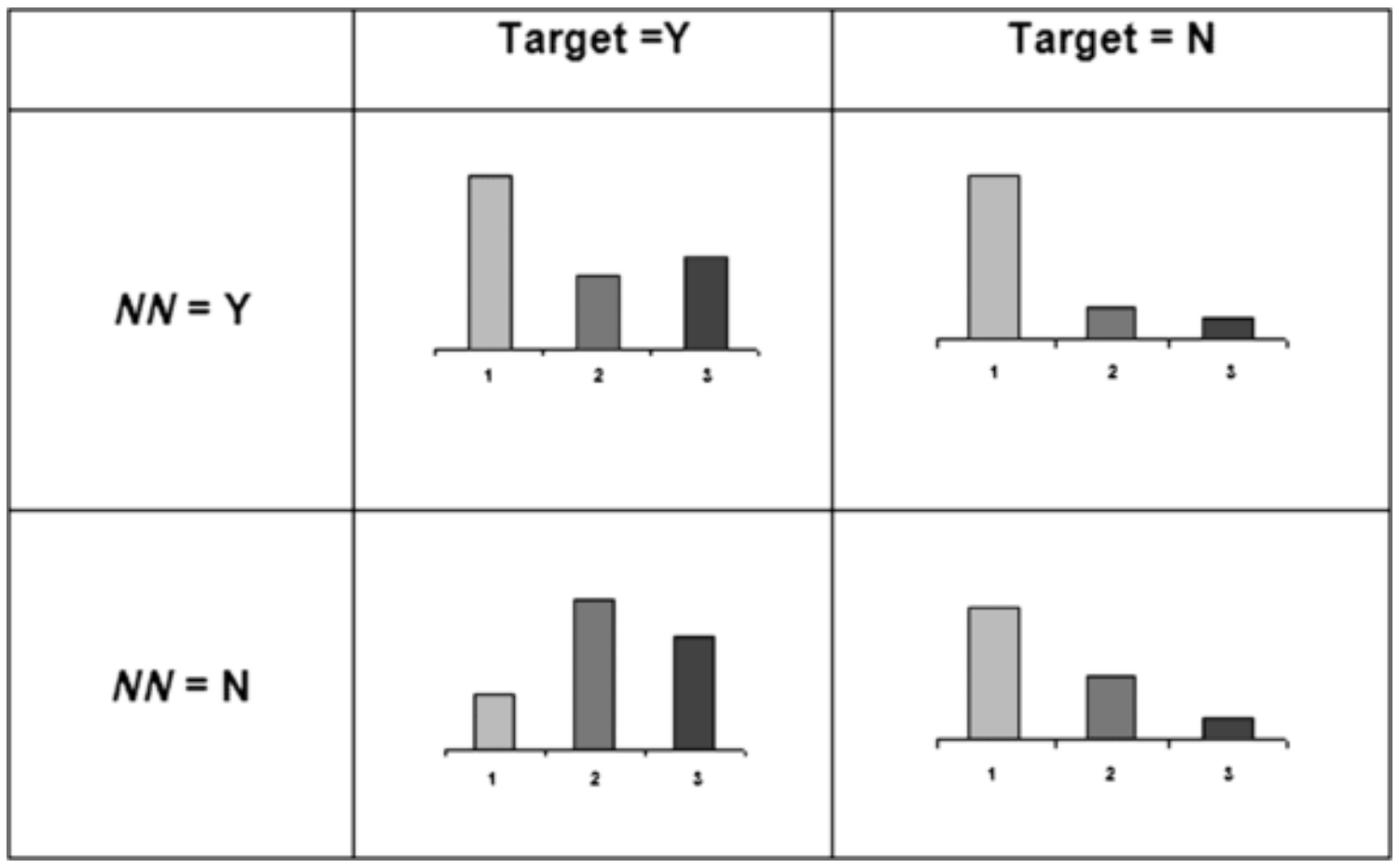

Fig. 1. *Loyalty distribution by customer category. The three bars show the relative sizes of the three loyalty groups for our four data categories. The network neighbors (NN) show a much larger proportion of low-loyalty consumers than the non-NN group.*

As Figure 3 shows, some of the segments had much higher take rates than others. To assess statistical significance of the network-neighbor effect after accounting for this segment effect, we ran a logistic regression across all segments, including the main effects for the network-neighbor attribute, dummy attributes for each segment and the interaction terms between the two. Two of the interaction terms had to be deleted: one from segment 22, which only had network-neighbor cases, and one from the segment with no sales from the network neighbors. We ran a full logistic regression and used stepwise variable selection.

The results of the logistic regression reiterate the significance of being a network neighbor. The final model can be found in Table 3. The coefficient of 2.0 for the network-neighbor attribute in the final model is an estimate of the log odds, which we exponentiate to get an odds ratio of 7.49, with a 95% confidence interval of (5.64, 9.94). More than half of the segment effects and most of the interactions between the network-neighbor attribute and those segment effects are significant. The interpretation of these interactions is important. Note that the magnitudes of the interaction coefficients are negative and very close in magnitude to the coefficients of the main effects of the segments themselves. Therefore, although the segments themselves are significant, in the presence of the network attribute the segments' effect is mostly negated by the interaction effect. Since the segments represent known important attributes like loyalty, tenure and demographics, this is evidence that being a network neighbor is at least as important in this context.

Table 3
*Coefficients and confidence intervals for the final segment model*

| **Attribute** | **Coeff (c.i.)** | **Significance**[a] |
|---|---|---|
| Network neighbor (NN) | 2.0 (1.7, 2.3) | ** |
| Segment = 1 | 1.7 (0.9, 2.5) | ** |
| Segment = 2 | 1.8 (1.2, 2.4) | ** |
| Segment = 4 | 2.1 (1.3, 3.0) | ** |
| Segment = 5 | 1.9 (0.4, 3.3) | ** |
| Segment = 6 | 1.9 (1.2, 2.5) | ** |
| Segment = 7 | 1.4 (1.0, 1.9) | ** |
| Segment = 8 | 1.3 (0.9, 1.7) | ** |
| Segment = 17 | 1.5 (0.7, 2.2) | ** |
| Segment = 19 | 2.2 (1.6, 2.9) | ** |
| NN × Segment = 1 | −1.1 (−2.1, 0.0) | * |
| NN × Segment = 2 | −0.9 (−1.7, −0.2) | ** |
| NN × Segment = 4 | −1.8 (−4.0, 0.4) | ** |
| NN × Segment = 6 | −1.5 (−2.6, −0.6) | ** |
| NN × Segment = 7 | −1.2 (−1.7, −0.6) | ** |
| NN × Segment = 8 | −0.8 (−1.3, −0.4) | ** |
| NN × Segment = 17 | −1.6 (−2.8, −0.5) | ** |
| NN × Segment = 19 | −1.1 (−1.9, −0.3) | ** |

[a] Significance of the attributes in the logistic regression model is shown at the 0.05 (*) and 0.01 (**) levels.

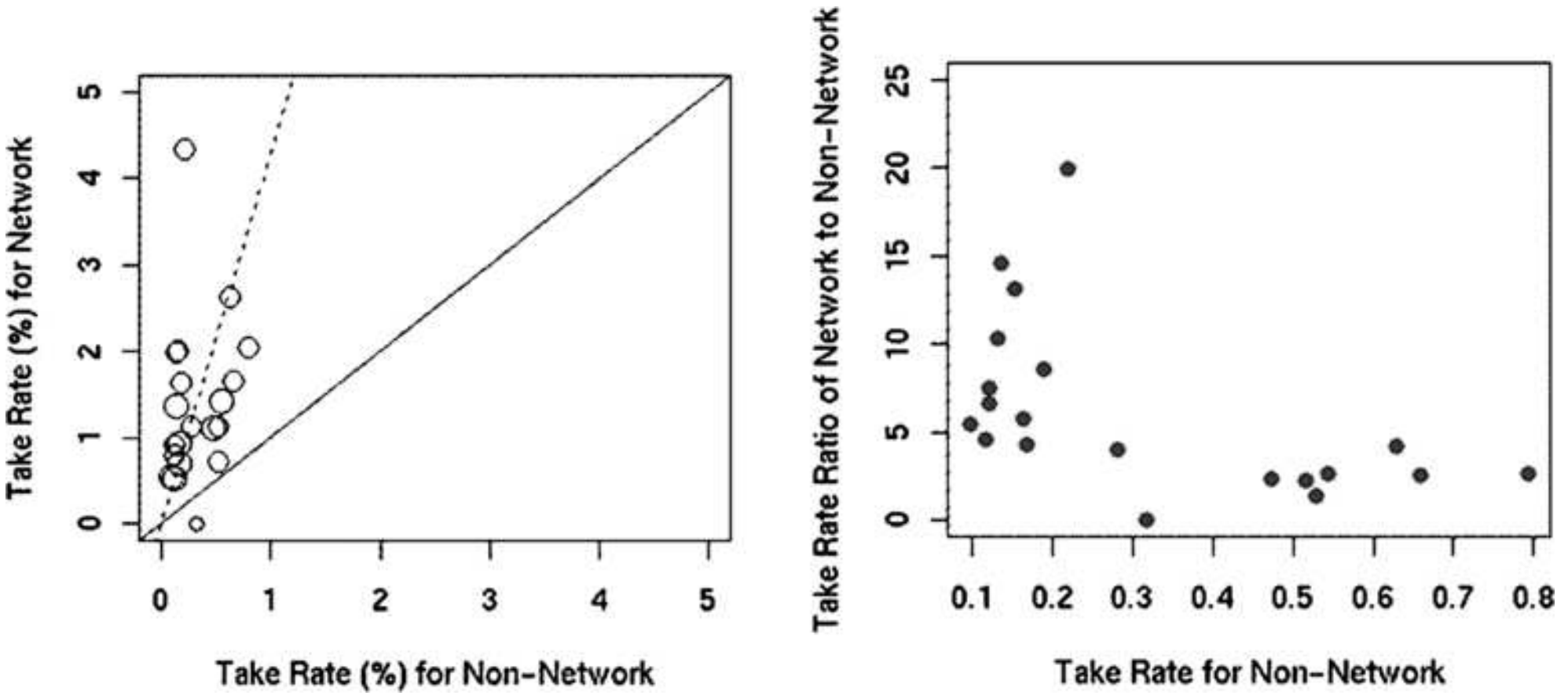

FIG. 3. *Take rates for marketing segments. Left: For each segment, comparison of the take rate of the non-network neighbors with that of the network neighbors. The size of the glyph is proportional to the log size of the segment. There is one outlier not plotted, with a take rate of* 11% *for the network neighbors and* 0.3% *for the non-network neighbors. Reference lines are plotted at* $x = y$ *and at the overall take-rate ratio of* 3.4. *Right: Plot of the take rate for the non-network group versus lift ratio for the network neighbors.*

In Table 4 we present an analysis of deviance table, an analog to analysis of variance used for nested logistic regressions (McCullagh and Nelder, 1983). The table confirms the significance of the main effects and of the interactions. Each level of the nested model is significant when a chi-squared approximation is used for the differences of the deviances. The fact that so many interactions are significant demonstrates that the network-neighbor effect varies for different segments of the prospect population.

## 5.2 Segment 22

The segment data enable us to compare take rates of network and non-network targets for the segments that contained both types of targets. However, many of the network-neighbor targets fall into the network-only segment 22. Segment 22 comprises prospects that the original marketing models deemed not to be good candidates for targeting. As we can see from the distribution in Figure 1, this segment for the most part contains consumers who had no prior relationship with the firm.

We compare the take rates for segment 22 with the take rates for the combined group, including all of segments 1–21, in the leftmost three bars of Figure 4. The network-neighbor segment 22 is (not surprisingly) not as successful as the NN groups in segments 1–21, since the targets in segments 1–21 were selected based on characteristics that made them favorable for marketing. Interestingly, we see that the segment 22 network neighbors outperform the non-NN targets from segments 1–21. These segment 22 network neighbors, identified primarily on the basis of their network activity, were more likely by almost 3 to 1 to purchase than the more "favorable" prospects who were not network neighbors. Since those in segment 22 either were not identified by marketing analysts or were deemed to be unworthy prospects, they represent customers who would have "fallen through the cracks" in the traditional marketing process.

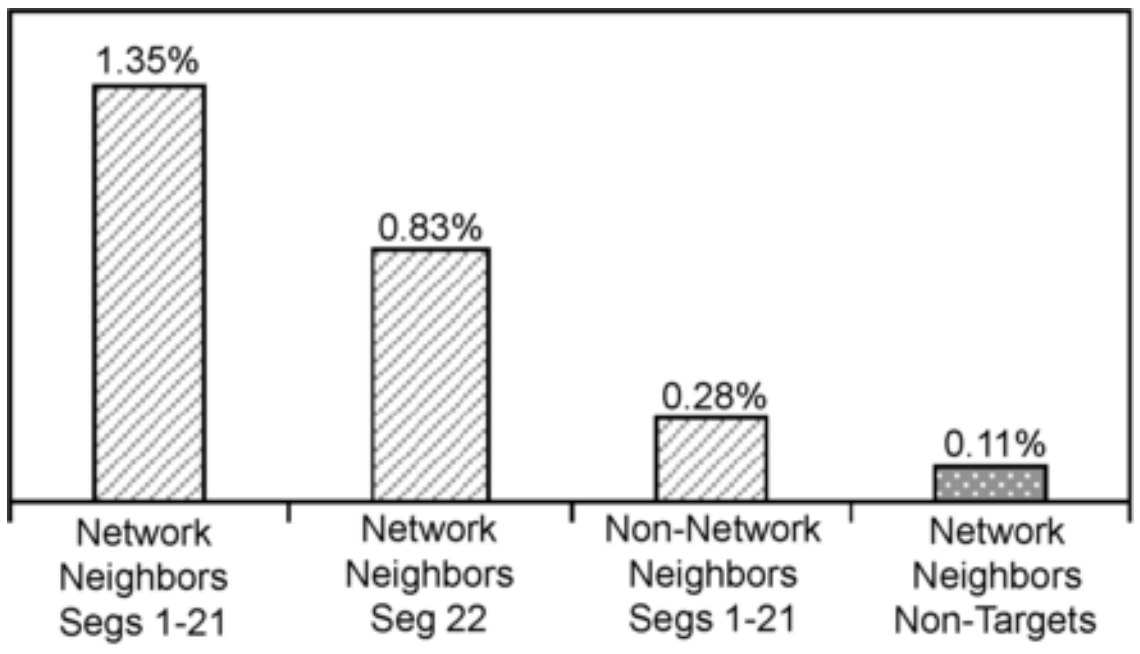

FIG. 4. *Take rates for marketing segments. Take rates for the network neighbors and non-network neighbors in segments* 1–21 *compared with the all-network-neighbor segment* 22 *and with the nontarget network neighbors. All take rates are relative to the non-network-neighbor group (segments* 1–21*).*

Table 4
*Analysis of deviance table for the network-neighbor study*

| Variable | Deviance | DF | Change in deviance | Significance[a] |
|---|---|---|---|---|
| Intercept | 11200 | | | |
| Segment | 10869 | 9 | 63 | ** |
| Segment + NN | 10733 | 1 | 370 | ** |
| Segment + NN + interactions | 10687 | 8 | 41 | ** |

[a]Significance of the group of attributes at each step is shown at the 0.05 (*) and 0.01 (**) levels.

## 5.3 Improving a Multivariate Targeting Model

Now we will assess whether the NN attribute can improve a multivariate targeting model by incorporating all that we know or can find out (over 150 different attributes) about the targets, including geography, demographics and other company-specific attributes, from internal and external sources (see Section 3.2).

As discussed in Section 3.7, we tried to address (as well as possible) an important causal question that arises: Is this network-neighbor effect due to word of mouth or simply due to homophily? The observed effect may not be indicating viral propagation, but instead may simply demonstrate a very effective way to find like-minded people. This theoretical distinction may not matter much to the firm for this particular type of marketing process, but is important to make, for example, before designing future campaigns that try to take advantage of word-of-mouth behavior.

Although we cannot control for unobserved similarities, we can be as careful as possible in our analysis to ensure that the statistical profile of the NN prospects is the same as the profile for the non-NN cases. Since our data set contains many more non-NN cases than NN cases, we match each NN case with a single non-NN case that is as close as possible to it by calculating propensity scores using all of the explanatory attributes considered (as described in Section 3.7). At the end of this matching process, the NN group is as close as is reasonably possible in statistical properties to the non-NN group.

Due to heterogeneity of data sources across the three loyalty groups, we used the propensity scores to create a matched data set for each group. For each (individually), we fitted a full logistic regression including interactions and selected a final model using stepwise variable selection. All attributes were checked for outliers, transformations and collinearity with other attributes, and we removed or combined the attributes that accounted for any significant correlations.

Table 5 shows the results of the logistic regressions, which show the attributes that were found to be significant, those that were negatively correlated with take rate, and those that had interactions with the NN attribute. Each of the three models found the network-neighbor attribute to be significant along with several others. The significant attributes tended to be attributes regarding the prospects' previous relationships with the firm, such as previous international services, tenure with firm, churn identifiers and revenue spent with the firm. These attributes are typically correlated with demographic attributes, which explains the lack of significance of many of the demographic attributes considered. Interestingly, tenure with firm is significant in loyalty groups 1 and 2, but with different signs. In the most loyal group, tenure is negatively correlated, but in the mid-level loyalty group it is positive. This unexpected result may be due to differing compositions of the two groups; those consumers with long tenure in the most loyal group might be people who just never change services, while long tenure in the other group might be an indicator that they are gaining more trust in the company. In loyalty group 1, there is limited information about previous services with the firm. For those customers, knowing whether the customer has responded to any previous marketing campaigns has a significant effect.

Table 5 also shows parameter estimates for NN and the take rates in the three loyalty groups. The take rates are highest in the group with the most

TABLE 5
*Results of multivariate model*

| | Loyalty | | |
|---|---|---|---|
| | **3** | **2** | **1** |
| **Significant attributes** | **NN** | **NN** | **NN** |
| | **Discount calling plan (-)(I)** | **Discount calling plan (-)** | **Previous responder to mailing** |
| | **Level of Int'l Comm.(I)** | **Tenure with firm** | High Tech Msg |
| | **# of devices in house (-)** | **Referral plan** | Letter (vs. postcard) |
| | **Revenue band** | **High Tech model score (I)** | Recent responder to mailing |
| | Tenure with firm (-) | Region of country indicator | User of incentive credit card |
| | International communicator | Belonged to loyalty program | Any children in house (-) |
| | Belonged to loyalty program | Churner (-) | |
| | Referral plan | College grad | |
| | Type of previous service | Tenure at residence (-) | |
| | Credit score | Any children in house (-) | |
| | Number of adults in house | Child < 18 at home (-) | |
| **Beta hat for NN (95% CI)** | **0.68 (0.46, 0.91)** | **0.99 (0.49, 1.49)** | **0.84 (0.52, 1.16)** |
| **Take rate** | **0.9%** | **0.4%** | **0.3%** |

NOTES. Significant attributes from logistic regressions across loyalty levels ($p < 0.05$). Bold indicates significance at 0.01 level; (-) indicates the effect of the variable was negative; (I) indicates a significant interaction with the NN variable.

loyalty but, interestingly, this group gets the least lift (smallest parameter estimate) from the NN attribute. So the impact of network-neighbor is stronger for those market segments with lower loyalty, where actual take rates are weakest.

## 5.4 Consumers Not Targeted

As discussed above, only a select subset of our network-neighbor list was subject to marketing, based on relaxed thresholds on eligibility criteria. The remainder of the list, the *nontarget network neighbors*, made up the majority. Potential customers were omitted for various reasons: they were not believed to have high-tech capacity; they were on a do-not-contact list; address information was unreliable, and so on. Nonetheless, we were able to identify whether they purchased the product in the follow-up time period. The take rate for this group was 0.11%, and is shown relative to the target groups as the rightmost bar in Figure 4. Although *they were not even marketed to,* their take rate is almost half that for the non-NN targets—chosen as some of the best prospects by the marketing team. This group comprises consumers without any known favorable characteristics that would have put them on the list of prospects. The fact that they are network neighbors alone supports a relatively high take rate, even in the absence of direct marketing. This lends some support to an explanation of word-of-mouth propagation rather than homophily.

Finally, we will briefly discuss the remainder of the consumer space—the non-NN nontarget group. Unfortunately, it is very difficult to estimate a take rate in this category, which could be considered a baseline rate for all of the other take rates. To do this, we would need to estimate the size of the space of all prospects. This includes all of the prospects the firm knows about, as well as customers of the firm's competitors and consumers who might purchase this product that do not have current telecommunications service with any provider. It has been established that the size of the communications market is difficult to estimate (Poole, 2004); our best estimates of this baseline take rate put it at well below 0.01%, at least an order of magnitude less than even the nontarget network neighbors.

On the other hand, a by-product of our study is that we can upper-bound the effect of the mass marketing campaigns in general by comparing the target-NN group and the nontarget-NN group. The difference in take rates between the targeted network neighbors and the nontargeted network neighbors is about 10 to 1. This difference cannot all be attributed to the marketing effect, since the targeted group was specifically chosen to be better prospects and it is likely that more of them would have signed up for the service even in the total absence of marketing. However, it does seem reasonable to call this factor of 10 an upper bound on the effect of the marketing.

## 5.5 Out-of-Sample Ranking Performance

These results suggest that we can give fine-grained estimations as to which customers are more or less likely to respond to an offer. Such estimations can be quite valuable: the consumer pool is immense and a campaign will have a limited budget. Therefore, being able to pick a better list of "top-$k$" prospects will lead directly to increased profit (assuming targeting costs are not much higher for higher ranked prospects). In this section, we show that combining the network-neighbor attribute with the traditional attributes improves the ability to rank customers accurately.

For each consumer, we create a record that comprises all of the traditional attributes (trad atts), including loyalty, demographic and geographic attributes, as well as network-neighbor status. Note that in different business scenarios, different types and amounts of data are available. For example, for low-loyalty customers, very few descriptive attributes are known. We report results here using all attributes; the findings are qualitatively similar for every different subset of attributes we have tried (namely, segment, loyalty, geography, demographic). The response variable is the same as above and we used the same logistic regression models. We measure the predictive ranking ability in the binary response variable by an increase in the Wilcoxon–Mann–Whitney statistic, equivalent to the area under the ROC curve (AUC). The ROC curve represents the trade off between false negative and false positive rates for each predicted possible probability score cutoff resulting from the logistic regression model. Specifically, the AUC is the probability that a randomly chosen (as yet unseen) taker will be ranked higher than a randomly chosen nontaker; AUC = 1.0 means the classes are perfectly separated and AUC = 0.5 means the list is randomly shuffled. All reported AUC values are averages obtained using 10-fold cross-validation.

Table 6 shows the AUC values for the three loyalty groups, quantifying the expected benefit from the improved logistic regression models. There is an increase in AUC for each group, with the largest increase belonging to loyalty level 1, for which the least information is available; note that here the ranking ability without the network information is not much better than random.

To visualize this improvement, Figure 5(a) shows cumulative response ("lift") curves when using the model on loyalty group 3. The lower curve depicts the performance of the model using all traditional attributes, and the upper curve includes the traditional marketing attributes and the network-neighbor attribute. In Figure 5(b), one can see the marked improvement that would be obtained from sending to the top-$k$ prospects on the list. For example, for the top 20% of the list, without the NN attribute, the take rate is 1.51%; with the NN attribute, it is 1.72%. The NN attribute does not improve the ranking for the top 10% of the list.

## 5.6 Improving Performance By Adding More Sophisticated Network Attributes

Knowing whether a consumer is a network neighbor is one of the simplest indicators of consumer-to-consumer interaction that can be extracted from the network data. We now investigate whether augmenting the model with more sophisticated social-network information can add additional value. In this section, we focus on the social network that comprises (only) the current customers of this service (which here we will call "the network"), along with the periphery of prospects who have communicated with those on the network (the network neighbors). We investigate whether we can improve targeting by using more sophisticated measures of social relationship with the network of existing customers.

Table 7 summarizes a set of additional social-network attributes that we add to the logistic regression. The terminology we use is borrowed to some degree from the fields of *social-network analysis* and *graph theory*. Social-network analysis (SNA) involves measuring relationships (including information transmission) between people on a network. The nodes in the network represent people and the links between

TABLE 6
*ROC analysis: AUC values that result from the application of logistic regression models*

| Loyalty | trad atts | trad atts + NN |
|---|---|---|
| 1 | 0.54 | 0.60 |
| 2 | 0.64 | 0.67 |
| 3 | 0.60 | 0.64 |

NOTE. The logistic regression models were built using all available attributes with (trad atts + NN) and without (trad atts) the network-neighbor attribute. We see an increase in AUC across all loyalty groups when the NN attribute is included in the model.

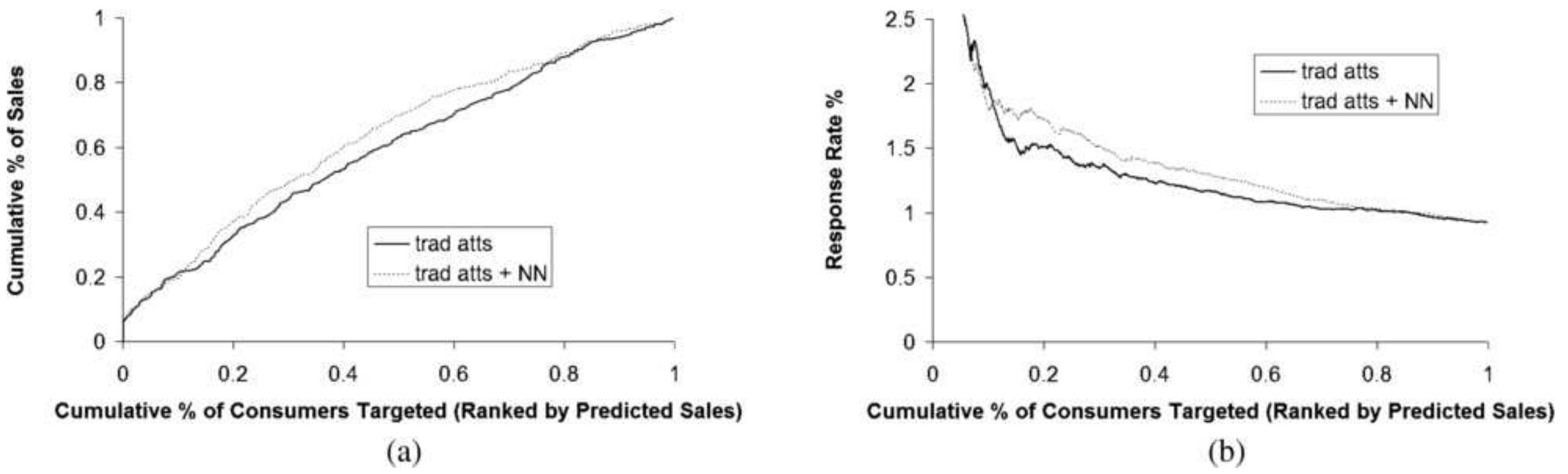

FIG. 5. *(a) Lift curves. Power of the segmentation curves for models built with all attributes with (trad atts) and without (trad atts + NN) network-neighbor attribute. The model with the NN attribute outperforms the model without it. For example, if the firm sent out* 50% *of the mailing, they would get* 70% *of the positive responses with the NN compared to receiving only* 63% *of the responses without it. (b) Top-k analysis. Consumers are ranked by the probability scores from the logistic regression model. The model that includes the NN attribute outperforms the model without. For example, for the top* 20% *of targets, the take rate is* 1.51% *without the NN attribute and* 1.72% *with the NN attribute.*

them represent relationships between the nodes. The SNA measures help quantify intuitive social notions, such as connectedness, influence, centrality, social importance and so on. Graph theory helps to understand problems better by representing them as interconnected nodes, and provides vocabulary and methods for operating mathematically.

Three of the attributes that we introduce can be derived from a prospect's local neighborhood (the set of immediate communication partners on the network; recall that these all are current customers). *Degree* measures the number of direct connections a node has. Within the local neighborhood, we also count the number of *Transactions*, and the length of those transactions (*Seconds of communication*).

The network is made up of many disjoint subgraphs. Given a graph $G = (V, E)$, where $V$ is a set of vertices (nodes) and $E$ is a set of links between them, the *connected components* of $G$ are the sets of vertices such that all vertices in each set are mutually connected (reachable by some path) and no two vertices in different sets are connected. The size of the connected component may be an indicator for awareness of and positive views about the product. If a prospect is linked to a large set of "friends" all of whom have adopted the service, she may be more likely to adopt herself. *Connected component size* is the size of the largest connected component (in the network) to which the prospect is connected.

We also move beyond a prospect's local neighborhood. Observing the local neighborhoods of a prospect's local neighbors, we can define a measure of social similarity. We define *social similarity* as the size of the overlap in the immediate network neighborhoods of two consumers. *Max similarity* is the maximum social similarity between the prospect and any neighbors of the prospect. Finally, the firm also can observe the prior dynamics of its customers. In particular, the firm can observe which customers communicated before and/or after their adoption as well as the date customers signed up. Using this information, we define *influencers* as those subscribers

TABLE 7
*Network attribute descriptions*

| **Attribute** | **Description** |
|---|---|
| Degree | Number of unique customers communicated with before the mailing |
| Transactions | Number of transactions to/from customers before the mailing |
| Seconds of communication | Number of seconds communicated with customers before mailing |
| Connected to influencer | Is an influencer in prospect's local neighborhood? |
| Connected component size | Size of the connected component prospect belongs to |
| Max similarity | Max overlap in local neighborhood with any existing neighboring customer |

TABLE 8
*ROC analysis*

| **Attribute(s)** | **AUC** |
|---|---|
| Transactions | 0.68 |
| Seconds of communication | 0.68 |
| Degree | 0.59 |
| Connected to influencer | 0.53 |
| Connected component size | 0.55 |
| Similarity | 0.55 |
| All network | 0.71 |
| All traditional (loyalty, demographic, geographic) | 0.66 |
| All traditional + all network | 0.71 |

NOTE. AUC values result from logistic regression models built on each of the constructed network attributes individually, as well as in combination. Results are presented for loyalty-level 3 customers.

who signed up for the service and, subsequently, we see one of their network neighbors sign up for the service. *Connected to influencer* is an indicator of whether the prospect is connected to one of these influencers. We appreciate that we do not actually know if there was true influence.

We use all of the aforementioned attributes and show AUC values for these predictive models in Table 8. We find that some of these network attributes have considerable predictive power individually and have even more value when combined. This is indicated by AUCs of 0.68 for both transactions and seconds of communication. We do not find high AUCs individually for connected component size, similarity or connected to influencer. Ultimately, we find that the logistic regression model built with the network attributes results in an AUC of 0.71 compared to an AUC of 0.66 without the network attributes—using only the traditional marketing attributes described in previous sections. (Recall that this represents the ability to rank *the network neighbors*, who already have especially high take rates as a group, as we have shown.)

Interestingly, when we combine the traditional attributes with the network attributes, there is no additional gain in AUC, even though many of these attributes were shown to be significant in the broader analysis above. The similarities represented implicitly or explicitly in the network attributes seem to account for all useful information captured by traditional demographics and other marketing attributes. That traditional demographics and other marketing attributes do not add value is not only of theoretical interest, but practical as well—for example, in cases such as this where demographic data must be purchased.

Our result is further confirmed by the lift and take rate curves displayed in Figure 6(a) and (b), respectively. One can achieve substantially higher take rates using the new network attributes as compared to using the traditional attributes. For example, we find that for the top 20% of the targeted list, without the network attributes, the take rate is 2.2%; with the network attributes, it is 3.1%. Likewise, at the top 10% of the list, the take rate with the network attributes is 4.4% compared to 2.9% without them.

## 6. LIMITATIONS

We believe our study to be the first to combine data on direct customer communication with data on product adoption to show the effect of network-based marketing statistically. However, there are limitations in our study that are important to point out.

There are several types of missing, incomplete or unreliable data which could influence our results. We have records of all of the communication (using the firm's service) to and from current customers of the service. That is not true for all the network-neighbor consumers. As such, we do not have complete information about the network-neighbor targets (as well as the non-network-neighbor targets). In addition, some of the attributes we used were collected by purchasing data from external sources. These data are known to be at least partially erroneous and outdated, although it is not well known how much so. An additional problem is joining data on customers from external sources to internal communication data, leading to missing data or sometimes just blatantly incorrect data. Finally, telecommunications firms are not legally able to collect information regarding the actual content of the communication, so we are not able to determine if the consumers in question discussed the product. In this regard, our data are inferior to some other domains where content is visible, such as Internet bulletin boards or product discussion forums.

We expect the network-neighbor effect to manifest itself differently for different types of products. Most of the studies done to date on viral marketing have focused on the types of products that people are likely to talk about, such as a new, high-tech gadget or a recently released movie. We expect there to be less buzz for less "sexy" products, like a new deodorant or a sale on grapes at the supermarket. The study presented in this paper involves a new telecommunications service, which involves a new technology and features that consumers have perhaps never been exposed to before. The firm hopes the new technology and features are such that they would encourage word of mouth.

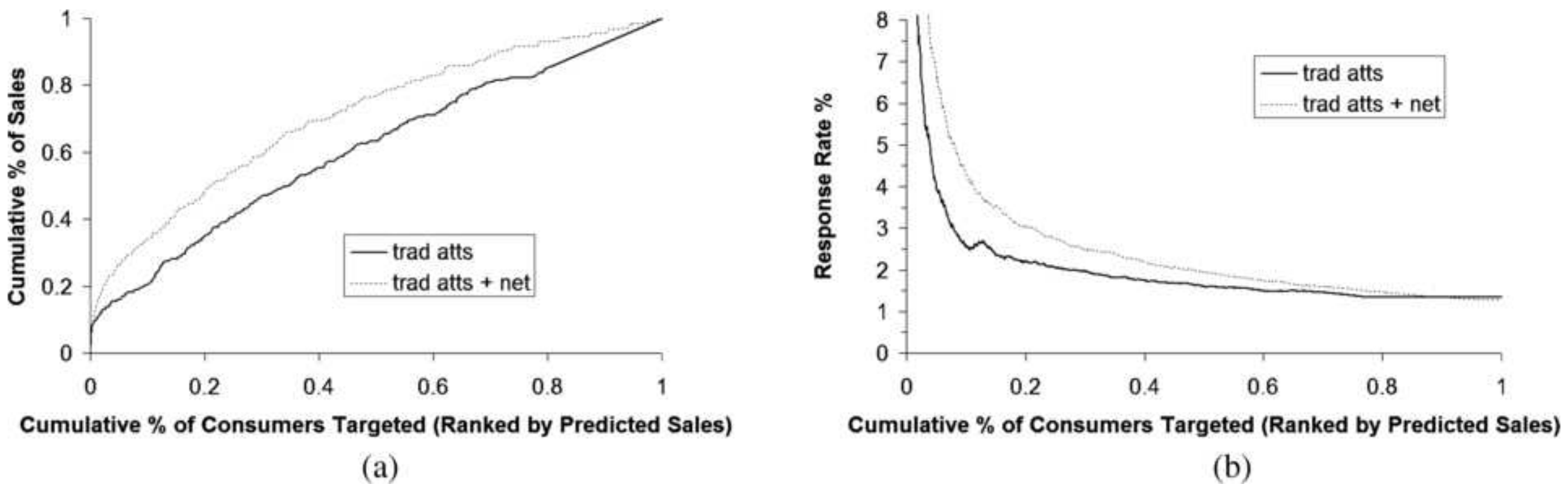

Fig. 6. (a) *Lift curves. Power of segmentation curves for models built with all traditional attributes, with (trad atts + net) and without (trad atts) the network attributes. If the firm sent out* 50% *of the mailing, they would have received* 77% *of the positive responses with the network attributes compared to receiving* 63% *of the responses without the network attributes.* (b) *Top-k analysis. The model including the network attributes (trad atts + net) outperforms the model without them (trad atts). For example, for the top* 20% *of target ranked by score, the take rate is* 2.2% *without the network attributes and* 3.1% *with the network attributes.*

What can we say about other products that might not be quite so buzz-worthy? To study this, we compared the new service studied here to a roll-out of another product by the same firm. This other product was simply a new pricing plan for an older telecommunications service. Customers who signed up for this new plan could stand to save a significant amount of money, depending on their current usage patterns. However, the range and variety of telecommunications pricing plans in the marketplace is so extensive and so confusing to the typical consumer that we do not believe that this is the type of product that would generate a lot of discussion between consumers. We refer to the two products as the *pricing plan* and the *new technology.* For the pricing plan, we have the same knowledge of the network as we do for the new technology. For those consumers who belong to the pricing plan, we know who they communicate with and then we can follow these network-neighbor candidates to see if they ultimately sign up for the plan. We construct a measure of "network neighborness" as follows. For a series of consecutive months, we gather data for all customers who ordered the product in that month. We calculate the percentage of these new customers who were network neighbors, that is, those who had previously communicated with a user of the product. This percentage is a measurement of the proportion of new sales being driven by network effects. By comparing this percentage across two products, we get insight into which product stimulates network effects more.

We now look at this value for our two products over an 8-month period. The time period for the two products was chosen so that it would be within the first year after the product was broadly available. The results are shown in Figure 7. The two main points to take away are that the new service has a higher percent of purchasers who are network neighbors and also an increasing one (except for the dip in month 5). In contrast the pricing plan has a flat network-neighbor percentage, never increasing above 3%.

Interestingly, the dip in the plot for the new service corresponds exactly to the month of the direct marketing discussed earlier. Before the campaign, we can see that the network-neighbor effect was increasing, that more and more of the purchasers in a given month were network neighbors. During the mass marketing campaign, we exposed many non-network neighbors to the service and many of them ended up purchasing it, temporarily dropping the network-neighbor percentage. After the campaign, we see the network-neighbor percentage starting to increase again.

This network-neighborness measure should not be confused with the success of the product, as the pricing plan was quite successful from a sales perspec-tive, but it does suggest that the pricing plan is a product that has less of a network-based spread of information. This difference might be due to the new service creating more word-of-mouth or perhaps we are seeing the effects of homophily. People who interact with each other are more likely to be similar in their propensity for purchasing the new service than in their propensity for purchasing a particular pricing plan. Again, the effects of word of mouth versus homophily are difficult to discern without knowing the content of the communication.

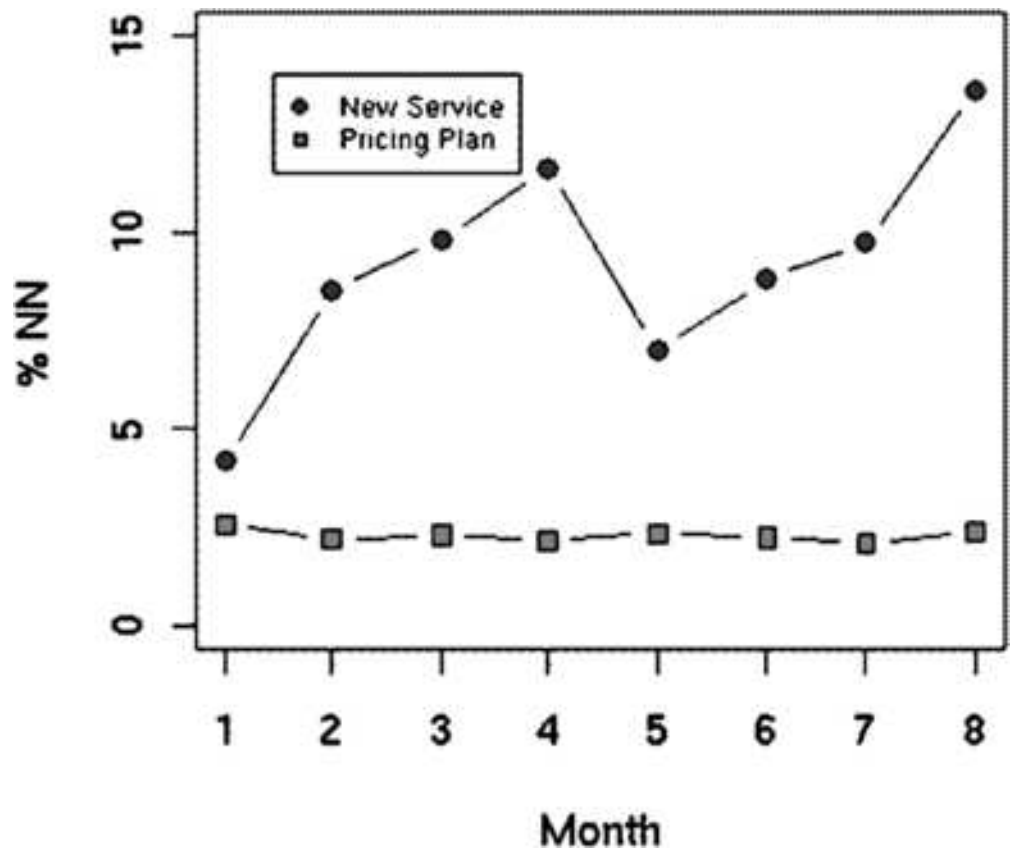

Fig. 7. *Network-neighborness plot for new service versus pricing plan.*

## 7. DISCUSSION

One of the main concerns for any firm is when, how and to whom they should market their products. Firms make marketing decisions based on how much they know about their customers and potential customers. They may choose to mass market when they do not know much. With more information, they may market directly based on some observed characteristics. We provide strong evidence that whether and how well a consumer is linked to existing customers is a powerful characteristic on which to base direct marketing decisions. Our results indicate that a firm can benefit from the use of social networks to predict the likelihood of purchasing. Taking the network data into account improves significantly and substantially on both the firm's own marketing "best practices" and our best efforts to collect and model with traditional data.

The sort of directed network-based marketing that we study here has applicability beyond traditional telecommunciations companies. For example, eBay recently purchased Internet-telephony upstart Skype for $2.6 billion; they now also will have large-scale, explicit data on who talks to whom. With gmail, Google's e-mail service, Google now has access to explicit networks of consumer interrelationships and already is using gmail for marketing; directed network-based marketing might be a next step. Various systems have emerged recently that provide explicit linkages between acquaintances (e.g., MySpace, Friendster, Facebook), which could be fruitful fields for network-based marketing. As more consumers create interlinked blogs, another data source arises. More generally, these results suggest that such linkage data potentially could be a sort of data considered for acquisition by many types of firms, as purchase data now are being collected routinely by many types of retail firms through loyalty cards. Even academic departments could benefit from such data; for example, the enrollment in specialized classes could be bolstered by "marketing" to those linked to existing students. Such links exist (e.g., via e-mail). It remains to design tactics for using them that are acceptable to all.

It is tempting to argue that we have shown that customers discuss the product and that discussion helps to improve take rates. However, word of mouth is not the only possible explanation for our result. As discussed in detail above, it may be that the network is a powerful source of information on consumer homophily, which is in accord with social theories (Blau, 1977; McPherson, Smith-Lovin and Cook, 2001). We have tried to control for homophily by using a propensity-matched sample to produce our logistic regression model. However, it may well be that direct communications between people is a better indicator of deep similarity than any demographic or geographic attributes. Either cause, homophily or word of mouth, is interesting both theoretically and practically.

## ACKNOWLEDGMENTS

We would like to thank DeDe Paul and Deepak Agarwal of AT&T, as well as Chris Dellarocas of the University of Maryland, for useful discussions and helpful suggestions. We would also like to thank three anonymous reviewers who offered insightful comments on previous drafts.

## REFERENCES

Adomavicius, G. and Tuzhilin, A. (2005). Toward the next generation of recommender systems: A survey of the state-of-the-art and possible extensions. *IEEE Trans. Knowledge and Data Engineering* **17** 734–749.

Agarwal, D. and Pregibon, D. (2004). Enhancing communities of interest using Bayesian stochastic blockmodels. In *Proc. Fourth SIAM International Conference on Data Mining*. SIAM, Philadelphia.

Bass, F. M. (1969). A new product growth for model consumer durables. *Management Sci.* **15** 215–227.

Blau, P. M. (1977). *Inequality and Heterogeneity*: *A Primitive Theory of Social Structure*. Free Press, New York.

Bowman, D. and Narayandas, D. (2001). Managing customer-initiated contacts with manufacturers: The impact on share of category requirements and word-of-mouth behavior. *J. Marketing Research* **38** 281–297.

Brin, S. and Page, L. (1998). The anatomy of a large-scale hypertextual Web search engine. *Computer Networks and ISDN Systems* **30** 107–117.

Case, A. C. (1991). Spatial patterns in household demand. *Econometrica* **59** 953–965.

Chan, E. and Stolfo, S. (1998). Toward scalable learning with non-uniform class and cost distributions: A case study in credit card fraud detection. In *Proc. Fourth International Conference on Knowledge Discovery and Data Mining* 164–168. AAAI Press, Menlo Park, CA.

Clearwater, S. H. and Stern, E. G. (1991). A rule-learning program in high-energy physics event classification. *Computer Physics Communications* **67** 159–182.

Dellarocas, C. (2003). The digitization of word of mouth: Promise and challenges of online feedback mechanisms. *Management Sci.* **49** 1407–1424.

Domingos, P. and Richardson, M. (2001). Mining the network value of customers. In *Proc. Seventh ACM SIGKDD International Conference on Knowledge Discovery and Data Mining* 57–66. ACM Press, New York.

DuMouchel, W., Volinsky, C., Johnson, T., Cortes, C. and Pregibon, D. (1999). Squashing flat files flatter. In *Proc. Fifth ACM SIGKDD International Conference on Knowledge Discovery and Data Mining* 6–15. ACM Press, New York.

Fichman, R. G. (2004). Going beyond the dominant paradigm for information technology innovation research: Emerging concepts and methods. *J. Assoc. Information Systems* **5** 314–355.

Fildes, R. (2003). Review of *New-Product Diffusion Models*, by V. Mahajan, E. Muller and Y. Wind, eds. *Internat. J. Forecasting* **19** 327–328.

Frenzen, J. and Nakamoto, K. (1993). Structure, cooperation, and the flow of market information. *J. Consumer Research* **20** 360–375.

Getoor, L. (2005). Tutorial on statistical relational learning. *Inductive Logic Programming*, *15th International Conference. Lecture Notes in Comput. Sci.* **3625** 415. Springer, Berlin.

Getoor, L., Friedman, N., Koller, D. and Pfeffer, A. (2001). Learning probabilistic relational models. In *Relational Data Mining* (S. Džeroski and N. Lavrač, eds.) 307–338. Springer, Berlin.

Getoor, L., Friedman, N., Koller, D. and Taskar, B. (2003). Learning probabilistic models of link structure. *J. Mach. Learn. Res.* **3** 679–707. MR1983942

Getoor, L. and Sahami, M. (1999). Using probabilistic relation models for collaborative filtering. In *Proc. WEBKDD 1999*, San Diego, CA.

Gladwell, M. (1997). The coolhunt. *The New Yorker* March 17, 78–88.

Gladwell, M. (2002). *The Tipping Point*: *How Little Things Can Make a Big Difference*. Back Bay Books, Boston.

Hightower, R., Brady, M. K. and Baker, T. L. (2002). Investigating the role of the physical environment in hedonic service consumption: An exploratory study of sporting events. *J. Business Research* **55** 697–707.

Hoff, P. D., Raftery, A. E. and Handcock, M. S. (2002). Latent space approaches to social network analysis. *J. Amer. Statist. Assoc.* **97** 1090–1098. MR1951262

Huang, Z., Chung, W. and Chen, H. C. (2004). A graph model for E-commerce recommender systems. *J. Amer. Soc. Information Science and Technology* **55** 259–274.

Japkowicz, N. and Stephen, S. (2002). The class imbalance problem: A systematic study. *Intelligent Data Analysis* **6** 429–449.

Joshi, M., Kumar, V. and Agarwal, R. (2001). Evaluating boosting algorithms to classify rare classes: Comparison and improvements. In *Proc. IEEE International Conference on Data Mining* 257–264. IEEE Press, Piscataway, NJ.

Kautz, H., Selman, B. and Shah, M. (1997). Referral web: Combining social networks and collaborative filtering. *Comm. ACM* **40**(3) 63–65.

Kleinberg, J. (1999). Authoritative sources in a hyperlinked environment. *J. ACM* **46** 604–632. MR1747649

Kumar, V. and Krishnan, T. V. (2002). Multinational diffusion models: An alternative framework. *Marketing Sci.* **21** 318–330.

Liben-Nowell, D. and Kleinberg, J. (2003). The link prediction problem for social networks. In *Proc. Twelfth International Conference on Information and Knowledge Management* 556–559. ACM Press, New York.

Linden, G., Smith, B. and York, J. (2003). Amazon.com recommendations—Item-to-item collaborative filtering. *IEEE Internet Computing* **7** 76–80.

Macskassy, S. and Provost, F. (2004). Classification in networked data: A toolkit and a univariate case study. CeDER Working Paper #CeDER-04-08, Stern School of Business, New York University.

Mahajan, V., Muller, E. and Kerin, R. (1984). Introduction strategy for new products with positive and negative word-of-mouth. *Management Sci.* **30** 1389–1404.

McCullagh, P. and Nelder, J. A. (1983). *Generalized Linear Models*. Chapman and Hall, New York. MR0727836

McPherson, M., Smith-Lovin, L. and Cook, J. (2001). Birds of a feather: Homophily in social networks. *Annual Review of Sociology* **27** 415–444.

Mease, D., Wyner, A. and Buja, A. (2006). Boosted classification trees and class probability/quantile estimation. *J. Mach. Learn. Res.* To appear.

Montgomery, A. L. (2001). Applying quantitative marketing techniques to the Internet. *Interfaces* **31**(2) 90–108.

Newton, J. and Greiner, R. (2004). Hierarchical probabilistic relational models for collaborative filtering. In *Proc.*

*Workshop on Statistical Relational Learning*, *21st International Conference on Machine Learning*. Banff, Alberta, Canada.

Paumgarten, N. (2003). No. 1 fan dept. acknowledged. *The New Yorker* May 5.

Perlich, C. and Provost, F. (2006). Distribution-based aggregation for relational learning with identifier attributes. *Machine Learning* **62** 65–105.

Poole, D. (2004). Estimating the size of the telephone universe. A Bayesian Mark-recapture approach. In *Proc. Tenth ACM SIGKDD International Conference on Knowledge Discovery and Data Mining* 659–664. ACM Press, New York.

Richardson, M. and Domingos, P. (2002). Mining knowledge-sharing sites for viral marketing. In *Proc. Eighth ACM SIGKDD International Conference on Knowledge Discovery and Data Mining* 61–70. ACM Press, New York.

Rogers, E. M. (2003). *Diffusion of Innovations*, 5th ed. Free Press, New York.

Rosenbaum, P. R. and Rubin, D. B. (1984). Reducing bias in observational studies using subclassification on the propensity score. *J. Amer. Statist. Assoc.* **79** 516–524.

Tout, K., Evans, D. J. and Yakan, A. (2005). Collaborative filtering: Special case in predictive analysis. *Internat. J. Computer Mathematics* **82** 1–11. MR2159280

Ueda, T. (1990). A study of a competitive Bass model which takes into account competition among firms. *J. Operations Research Society of Japan* **33** 319–334.

Van den Bulte, C. and Lilien, G. L. (2001). Medical innovation revisited: Social contagion versus marketing effort. *American J. Sociology* **106** 1409–1435.

Walker, R. (2004). The hidden (in plain sight) persuaders. *The New York Times Magazine* Dec. 5, 69–75.

Weiss, G. and Provost, F. (2003). Learning when training data are costly: The effect of class distribution on tree induction. *J. Artificial Intelligence Research* **19** 315–354.

Weiss, G. M. (2004). Mining with rarity: A unifying framework. *ACM SIGKDD Explorations Newsletter* **6** 7–19.

Yang, S. and Allenby, G. M. (2003). Modeling interdependent consumer preferences. *J. Marketing Research* **40** 282–294.